\documentclass[letterpaper, 10 pt, conference]{ieeeconf}  

\IEEEoverridecommandlockouts                              

\usepackage{graphicx}      
\usepackage{amsmath}
\usepackage[utf8]{inputenc}
\usepackage{graphicx}
\usepackage{algorithmicx}
\usepackage{algorithm}
\usepackage{algpseudocode}
\usepackage{algorithm,algpseudocode}
\usepackage{amssymb}
\usepackage{float}
\usepackage{enumerate}
\usepackage{comment}
\usepackage{subfigure}
\usepackage[table,xcdraw]{xcolor}

\newtheorem{Problem}{\textbf{Problem}}

\title{\LARGE \bf
Quadratic approximation based heuristic for optimization-based coordination of automated vehicles in confined areas
}

\author{Stefan Kojchev$^{1}$, Robert Hult$^{2}$ and Jonas Fredriksson$^{3}$
\thanks{*This work is partially funded by Sweden's innovation agency Vinnova, project number: 2018-02708.}
\thanks{$^{1}$Stefan Kojchev is with Volvo Autonomous Solutions and the Mechatronics Group, Systems and Control, Chalmers University of Technology
        {\tt\small stefan.kojchev@volvo.com}; {\tt\small kojchev@chalmers.se}}%
\thanks{$^{2}$Robert Hult is with Volvo Autonomous Solutions, 41873 Göteborg, Sweden
        {\tt\small robert.hult@volvo.com}}%
\thanks{$^{3}$Jonas Fredriksson is with the Mechatronics Group, Systems and Control, Chalmers University of Technology, 41296 Göteborg, Sweden
        {\tt\small jonas.fredriksson@chalmers.se}}%
}

\newcommand{\tder}[2]{\frac{\mathrm{d}#1}{\mathrm{d}#2}}

\begin{document}


\maketitle
\thispagestyle{empty}
\pagestyle{empty}

\begin{abstract}
We investigate the problem of coordinating multiple automated vehicles (AVs) in confined areas. This problem can be formulated as an optimal control problem (OCP) where the motion of the AVs is optimized such that collisions are avoided in cross-intersections, merge crossings, and narrow roads. The problem is combinatorial and solving it to optimality is prohibitively difficult for all but trivial instances. For this reason, we propose a heuristic method to obtain approximate solutions. The heuristic comprises two stages: In the first stage, a Mixed Integer Quadratic Program (MIQP), similar in construction to the Quadratic Programming (QP) sub-problems in Sequential Quadratic Programming (SQP), is solved for the combinatorial part of the solution. In the second stage, the combinatorial part of the solution is held fixed, and the optimal state and control trajectories for the vehicles are obtained by solving a Nonlinear Program (NLP). The performance of the algorithm is demonstrated by a simulation of a non-trivial problem instance. 
\end{abstract}

\section{Introduction}\label{Ch:Introduction}

The idea of fully automated vehicles (AV) is receiving substantial attention from both the public and scientific world as significant progress towards deploying automated vehicles has been made during the last decade. Unfortunately, many barriers between the current state-of-the-art and large-scale commercial application exits, especially for deployment of automated vehicles on public roads \cite{b1}. However, in confined areas, such as ports, mines, and logistic centers, some of the hard challenges of public road driving are absent. In particular, such areas are typically void of unpredictable non-controlled actors, which dramatically reduces safety concerns. Therefore, it is believed that confined areas present a good opportunity for early, large-scale deployment of automated vehicles, as part of larger commercial transport solutions for material flow. 

One of the challenges in confined areas is the safe and efficient coordination of AVs in mutually exclusive (MUTEX) zones, such as intersections, work-stations (e.g. crushers, loading/unloading spots, etc.), narrow roads, and, in the case of electrified AVs, charging-stations. Adequate coordination can lead to improved energy-efficiency and considerable increases in productivity.

\subsection{Related Work}

Automating and coordinating intersections for fully automated vehicles is a frequently discussed control problem, see \cite{b2} for a comprehensive survey. The problem has been formally shown to be NP-hard \cite{b3}, and such problems are in general difficult to solve. For this reason a number of methods have been proposed, leveraging results from, e.g., solutions based on hybrid system theory \cite{b4}, reinforcement learning \cite{b5}, scheduling \cite{b6}, model predictive control (MPC) \cite{b7}, \cite{b8} or direct optimal control (DOC) \cite{b9}, \cite{b10}.

In contrast to intersection scenarios often found in the literature, confined areas have a number of distinguishing features. For example, in the case of intersection coordination, an approach that is often considered is to have vehicles arriving at speed in a cutout around the intersection area \cite{b7}, \cite{b8}. This is often motivated by practical considerations; neither the intent of the automated vehicles nor the state of the uncertain environment can be accurately predicted over long time-horizons. For confined areas, however, it is possible, and desirable, to plan the motion of each vehicle from the start of a transport mission to its end. Moreover, a number of works on public-road applications focus on distributed and decentralized schemes, sometimes with intermittent or corrupted communication \cite{b17}. For applications at confined sites, a central computational unit and good wireless coverage can often be assumed. Thus, for these use cases, we believe that a centralized approach that provides high level control actions is favorable. A low level controller that tracks the obtained optimal state and control trajectories will typically be also deployed in practice, however, it is not covered in this work. 
\begin{figure*}[htp]
    \centering
    \includegraphics[width=0.9\textwidth]{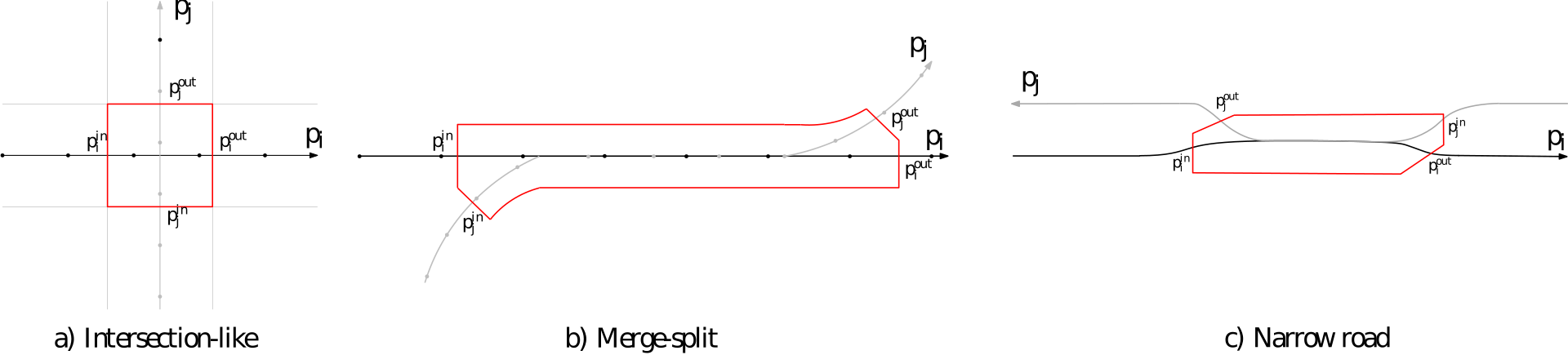}
    \caption{Types of conflict zones.}
    \label{fig:IntersectionsTogether}
\end{figure*}
\subsection{Main Contribution and Outline}

In this paper, we formulate the site-coordination problem as an optimal control problem, which after transcription results in a Mixed Integer Nonlinear Program, and propose a two-staged heuristic method for its approximate solution. In the first stage of the heuristic, an MIQP problem is solved to obtain the combinatorial part of the solution, i.e., the order in which the vehicles utilize the MUTEX-zones. In the second stage, the combinatorial part is fixed and a continuous NLP is solved for the optimal vehicle trajectories. The MIQP of the first stage is formed as an approximation of the original MINLP, in a manner similar to how the Quadratic Programming (QP) sub-problems are formed in Sequential Quadratic Programming (SQP) \cite{b11}. The idea of using SQP-like methods in approximate solution of MINLPs has been used in other works \cite{b13}, \cite{b14}, however, to the authors' best knowledge, it has not been adapted and applied to vehicle coordination problems. A structurally similar heuristic was presented in \cite{b9}, where a scheduling problem is derived from the original MINLP, with the introduction of a number of approximations, and solved as an MIQP.  While sharing the same two-staged structure, the method presented in this paper avoids some of the approximations and shortcomings of \cite{b9}, without expanding the combinatorial solution space. In particular,  the heuristic presented herein enables easy inclusion of rear-end collision constraints, which were previously neglected, and avoids the expensive computation of parametric sensitivities.

In addition to the cross-intersections, we consider merge-split and narrow road MUTEX zones. In the merge-split MUTEX zones, the vehicles first join in on a common patch of road which after some distance separate, and in the narrow-roads the vehicles that are coming from different directions join in on a common patch of road. Merge-split and narrow roads are often found in confined sites and occur due to the construction of the site. Although the approach focuses on confined sites, the method of handling the mutual exclusion zones can be extendable to other scenarios as well (e.g., public road applications). 

The remainder of the paper is organized as follows: Section \ref{Ch:Problem formulation} presents a formulation of the problem that is solved in this paper. In Section \ref{Ch:Method} the method for solving the stated problem is presented, followed by Section \ref{Ch:Simulation results} where simulation results illustrate the coordination algorithm. Section \ref{Ch:Conclusions} concludes the work and provides some possible extensions.

\section{Problem formulation}\label{Ch:Problem formulation}

In this paper, we consider a fully confined area, meaning that non-controlled traffic participants such as pedestrians, manually operated vehicles, bicycles, etc., are absent. Furthermore, the confined road network consists of $N_a$ fully automated vehicles with cross-intersection, path merges, path splits, and narrow roads. In addition, we assume that the paths of all vehicles, i.e., their routes through the road network are known, that overtakes are prohibited, and that no vehicle reverses.
\subsection{Optimal Coordination Problem}

The problem of finding the optimal vehicle trajectories that avoid collisions can be stated as: 

\Problem (Optimal coordination problem) Obtain the optimal state and control trajectories $\mathcal{X}^{*} = \left \{x_1^*,...,x_{N_a}^*\right \}$, $\mathcal{U}^{*} = \left \{u_1^*,...,u_{N_a}^*\right \}$, given the initial state $\mathcal{X}_{0} = \left \{x_{1,0},...,x_{N_a,0}\right \}$, by solving the optimization problem 
\begin{subequations}
\label{Eq: Main_OCP}
\begin{align} 
        \underset{x_{i},u_{i}, \mathcal{O}^\mathcal{I}, \mathcal{O}^\mathcal{M}}{\textup{min}} \;\; &\sum_{i = 1}^{N_a}J_i\left ( x_{i},u_{i} \right )  \\
\textup{s.t} \;\;\; & \textup{initial states} \; \; x_{i,0} = \hat{x}_{i,0}, \forall i \label{Eq: InitialGuess}  \\
&\textup{system dynamics}\; \; \forall i \label{Eq: OCPSysDyn} \\
&\textup{state and input constraints}  \; \; \forall i\label{Eq: OCPStateInpConst} \\
& \textup{safety constraints} \; \; \forall i  \label{Eq: OCPSafetyCosnt}
\end{align}
\end{subequations}
where $N_a$ is the number of vehicles, $J_i(x_i,u_i)$ is the cost function, $\mathcal{O}^\mathcal{I}, \mathcal{O}^\mathcal{M}$ are the order in which the vehicles enter the MUTEX zones and will be formally stated in this section.

The problem is formulated in the spatial domain as it is beneficial to optimize the trajectories of the vehicles over their full paths. The rationale for using spatial dynamics is that the time it takes for the vehicle to traverse a path is not known \textit{a-priori}. Thus, it is inappropriate to plan the vehicle's motion with time as the independent variable.

\subsection{System dynamics and state and input constraints}
The system dynamics for vehicle $i \in {1,...,N_a}$ in the spatial domain can be formed using that $\frac{\mathrm{d}p_i}{\mathrm{d}t}= v_i(t)$ and $\mathrm{d}t = \mathrm{d}p_i/v_i(t)$ and stated as
\begin{align}
    \tder{t_i}{p_i}&= \frac{1}{v_i(p_i)} \label{Eq: TimeDyn}\\
    \tder{x_i}{p_i} &= \frac{1}{v_i(p_i)}f_i(p_i,x_i(p_i),u_i(p_i)) \label{Eq: SpaceDyn}\\
    0 &\leq h(p_i,x_i,u_i) \label{Eq: State_Input_Const}.
\end{align}
\noindent where the position $p_i$ is the independent variable, the time $t_i$ and $x_i$ are the vehicle state variables, where $x_i \in \mathbb{R}^{n-1}$ collects the remaining vehicle states, and  $u_i\in\mathbb{R}^m$ the control input, with $i\in \{1, \hdots, N_a \}$. Note that what the remaining state variable ($x_i$) are, depends on what model is used for the system dynamics. We assume that the functions $f_i$ and $h_i$, that describe the vehicle system's dynamics and constraints, are smooth. 

\subsection{Safety constraints}

The safety constraints should ensure a collision-free crossing of the \textit{conflict zone} (CZ) that the vehicles encounter. A conflict zone is described by the entry and exit position $[p_i^{in},p_i^{out}]$ on the path of each vehicle. From the known positions, the time of entry and exit of vehicle $i$ is $t_i^{in}=t_i(p_i^{in})$ and $t_i^{out}=t_i(p_i^{out})$, respectively. In this paper, we consider three types of conflict zones, the ``intersection-like", ``narrow road" and the ``merge-split", depicted in Figure \ref{fig:IntersectionsTogether}. The narrow road CZ is when two vehicles that are coming from opposite directions have to share a common patch of road. In the intersection-like CZ and narrow road CZ, it is necessary to only have one vehicle inside the CZ, i.e., not allowing the vehicle $j$ to enter the CZ before vehicle $i \neq j$ exits the CZ, or vice-versa. We let $\mathcal{I} = \left \{ I_1, I_2,...,I_{r_0} \right \}$ denote the set of all intersections and narrow roads in the confined site, with $r_0$ being the total number of intersection and narrow road CZs, and $\mathcal{Q}_r = \left \{ q_{r,1}, q_{r,2},...,q_{r,l} \right \}$ denote the set of vehicles that cross an intersection or narrow road $I_r$. The order in which the vehicles cross the intersection $I_r$ is denoted $\mathcal{O}^{\mathcal{I}}_r = \left ( s_{r,1}, s_{r,2},...,s_{r,|\mathcal{Q}_r|} \right )$, where $s_{r,1}, s_{r,2},...$ are vehicle indices and we let $\mathcal{O}^\mathcal{I}=\left \{\mathcal{O}^\mathcal{I}_1, \hdots, \mathcal{O}^\mathcal{I}_r\right \}$. A sufficient condition for collision avoidance for the $r$-th intersection or narrow road CZ can be formulated as
\begin{align}\label{Eq: IntersectionConst}
    t_{s_{r,i}}(p_{s_{r,i}}^{out})\leq t_{s_{r,i+1}}(p_{s_{r,i+1}}^{in}), \; &i \in \mathbb{I}_{\left [ 1, |\mathcal{Q}_r|-1 \right ]},
\end{align}
\noindent where $t$ is determined from \eqref{Eq: TimeDyn}.

In the merge-split CZ case, let $\mathcal{M} = \left \{ M_1, M_2,...,M_{w_0}  \right \}$ denote a set of all merge-split zones, with $w_0 $ being the total number of merge-split CZs in the site and $\mathcal{Z}_w = \left \{ z_{w,1}, z_{w,2},...,z_{w,h} \right \}$ denote the set of vehicles that cross the merge-split CZ $M_w$. For efficiency, it is desirable to have several vehicles in the zone at the same time, instead of blocking the whole zone. This requires having rear-end collision constraints once the vehicles have entered the merge-split CZ. In this case, the order in which the vehicles enter the zone is denoted as $\mathcal{O}^{\mathcal{M}}_w = \left ( s_{w,1}, s_{w,2},...,s_{w,|\mathcal{Z}_w|}\right )$, and we let $\mathcal{O}^\mathcal{M}=\left \{\mathcal{O}^\mathcal{M}_1, \hdots, \mathcal{O}^\mathcal{M}_w\right \}$. The collision avoidance requirement for the $w$-th merge-split CZ is described with the following constraints: 
\begin{subequations}
\label{Eq:MergeSplit_ORConst}
\begin{align}\label{Eq: MS_In}
t_{s_{w,i}}(p_{s_{w,i}}^{in}) + \Delta t &\leq t_{s_{w,i+1}}(p_{s_{w,i+1}}^{in}+c) \\ \nonumber
t_{{s_{w,i}},k_i} + \Delta t &\leq t_{s_{w,i+1}}(p_{{s_{w,i}},k_i} - p_{s_{w,i}}^{in} + p_{s_{w,i+1}}^{in} + c), \\  & \;\;\;\;k_{s_{w,i}}^{in} \leq k_i \leq k_{s_{w,i}}^{out} \label{Eq: MS_Middle}  \\ 
t_{s_{w,i}}(p_{s_{w,i}}^{out}) + \Delta t &\leq t_{s_{w,i+1}}(p_{s_{w,i+1}}^{out}+c), \label{Eq: MS_Out}\\ \nonumber
&i \in \mathbb{I}_{\left [ 1, |\mathcal{Z}_w|-1 \right ]}.
\end{align}
\end{subequations}
That is, while in the CZ, the vehicles must be separated by at least a time-period $\Delta t$ and a distance $c$, depending on if vehicle $j$ is in front of vehicle $i$ or vice versa. This is equivalent to the standard offset and time-headway formulation often used in automotive adaptive cruise controllers.
\subsection{A practical reformulation of the collision constraints}
A common way to handle constraints such as \eqref{Eq: IntersectionConst} and \eqref{Eq:MergeSplit_ORConst} is to introduce auxiliary binary variables and use the ``big-M" technique \cite{b12}. For example, an equivalent representation to the constraint \eqref{Eq: IntersectionConst}, with $b_{s_r,i,i+1} \in \left\{ 0,1\right\}, \; i\in \mathbb{I}_{\left [ 1, |\mathcal{Q}_r|-1 \right ]}$ and a sufficiently large $M$ is
\begin{subequations}
\label{Eq: Big-M reformulation}
\begin{align}
     &t_{s_{r,i}}(p_{s_{r,i}}^{out}) - t_{s_{r,i+1}}(p_{s_{r,i+1}}^{in}) \leq b_{s_r,i,i+1}M, \\
     &t_{s_{r,i+1}}(p_{s_{r,i+1}}^{out}) - t_{s_{r,i}}(p_{s_{r,i}}^{in}) \leq (1-b_{s_r,i,i+1})M.
\end{align}
\end{subequations}
In the case where $b_{s_r,i,i+1}=0$, the vehicle $i+1$ is constrained to cross the CZ after the vehicle $i$, with the opposite being true if $b_{s_r,i,i+1}=1$. We collect all integer variables for all CZs in $b \in \mathbb{Z}_2^{r_o+w_0}$.
\subsection{Discretization}
The independent variable is discretized as $p_i = (p_{i,0}, \hdots, p_{i,N_i})$, where $p_{i,N_i}$ indicates the end position for vehicle $i$, and the input is approximated using zero order hold such that $u(p) = u_{i,k}, p\in[p_{i,k}, p_{i,k+1}[$. The equations (\ref{Eq: TimeDyn}), (\ref{Eq: SpaceDyn}) are (numerically) integrated on this grid, giving the ``discretized" state transition relation 
\begin{equation}\label{Eq: SystemDynamics}
    \begin{bmatrix}
    t_{i,k+1} \\
    x_{i,k+1}
    \end{bmatrix} = F(x_{i,k},u_{i,k},p_{i,k},p_{i,k+1}) 
\end{equation}
where $F$ denotes the integration of (\ref{Eq: TimeDyn}), (\ref{Eq: SpaceDyn}) from $p_{i,k}$ to $p_{i,k+1}$.
\section{Method}\label{Ch:Method}
The optimal coordination problem, Problem 1, can be stated as a Mixed Integer Nonlinear Program (MINLP), where the crossing order correspond to the combinatorial (``integer part") and the state and control trajectories correspond to the ``NLP part". In essence, we can state Problem 1 as
\begin{subequations}
\label{Eq: MINLP}
\begin{align}
&\underset{\mathcal{W},b}{\textup{min}} \;\; J(\mathcal{W}) \\ 
&\textup{s.t.} \;\;\; g(\mathcal{W}) = 0 \\
& \hspace{0.7cm} h(\mathcal{W})\leq 0 \\
& \hspace{0.7cm} c(\mathcal{W},b)\leq 0,
\end{align}
\end{subequations}
\noindent where $\mathcal{W} = \left \{\mathcal{X}, \mathcal{U}\right \}$, $J(\mathcal{W}) = \sum_{i = 1}^{N_a}J_i\left ( w_{i}\right )$, $g(\mathcal{W}),h(\mathcal{W})$ gather all equality and inequality constraints, and $c(\mathcal{W},b) = c_w(\mathcal{W}) + Cb$ are the integer constraints for the combinatorial part of the problem with $C$ being a matrix that captures the influence of the integer variables. 

Since finding a solution to MINLP problems is known to be difficult, it is common to obtain approximate solutions with heuristics. One heuristic approach used in e.g. \cite{b9}, is to use a two-staged procedure, where the integer part of the solution first is found with a heuristic, and the continuous part is found by solving the NLP obtained by eliminating all other integer options from the MINLP. In this paper, we follow this decomposition idea and propose an alternative heuristic for the integer part than \cite{b9}. In the confined area coordination context, the integer part of the solution corresponds to the crossing order of the vehicles at the MUTEX zones.

\subsection{Crossing order heuristic}
The crossing order heuristic is based on solving an MIQP that is assembled from a quadratic approximation of \eqref{Eq: MINLP}. The quadratic approximation is formed similarly to how the QP sub-problems are formed in SQP methods \cite{b11}. In essence, \eqref{Eq: MINLP} can be reformulated as: 
\begin{subequations}
\label{Eq: MIQP_SQP}
\begin{align} \nonumber
&\underset{\Delta \mathcal{W}, b}{\textup{min}} \;\; \frac{1}{2}\begin{bmatrix}
\Delta \mathcal{W} \\
b
\end{bmatrix}^T \mathbf{H}(\mathcal{W},\lambda,\mu) \begin{bmatrix}
\Delta \mathcal{W} \\
b
\end{bmatrix} + \\ 
&  \hspace{0.95cm} \triangledown_{\mathcal{W}} J(\mathcal{W})^T\begin{bmatrix}
\Delta \mathcal{W} \\
b
\end{bmatrix} + J(\mathcal{W}^{**})\\
&\textup{s.t.} \;\;\;\;\;\; g(\mathcal{W}^{**}) + \triangledown_{\mathcal{W}} g(\mathcal{W}^{**})^T\begin{bmatrix}
\Delta \mathcal{W}\\
b
\end{bmatrix} = 0 \label{Eq: SQP_Const1} \\
& \hspace{0.95cm} h(\mathcal{W}^{**})+ \triangledown_{\mathcal{W}} h(\mathcal{W}^{**})^T\begin{bmatrix}
\Delta \mathcal{W}\\
b
\end{bmatrix}\leq 0 \label{Eq: SQP_Const2}\\
& \hspace{0.95cm} c_w(\mathcal{W}^{**}) + \triangledown_{\mathcal{W}} c_w(\mathcal{W}^{**})^T\begin{bmatrix}
\Delta \mathcal{W}\\
b
\end{bmatrix} + Cb\leq 0, \label{Eq: SQP_Const3}
\end{align}
\end{subequations}
\noindent where $\mathbf{H}(\mathcal{W},\lambda,\mu) = \textup{blkdiag} (\left\{H_i\right\}_{i=1}^{N_a},\mathbf{0}_{r_0+w_0,r_0+w_0})$ is a block diagonal matrix with positive definite $H_i(w_i,\lambda_i,\mu_i) = \triangledown_{w_i}^2 \mathcal{L}(w_i,\lambda_i,\mu_i) = \triangledown_{w_i}^2J_i(w_i) - \triangledown_{w_i}^2\lambda_i^T g(w_i) - \triangledown_{w_i}^2\mu_i^T h(w_i)$, where $\lambda_i,\mu_i$ are the dual variables and $\mathbf{0}_{r_0+w_0,r_0+w_0}$ zeros of appropriate size for the integer variables, and $\Delta \mathcal{W} = \mathcal{W} - \mathcal{W}^{**}$, with a solution guess $\mathcal{W}^{**}$. 

For the heuristic used in this paper, we make the simplification that the dual variables ($\lambda_i,\mu_i)$ are equal to zero. This results in that $H_i$ only includes the second order expansion of the cost function, i.e., $H_i(w_i) = \triangledown^2_{w_i}J_i(w_i)$. The solution guess $\mathcal{W}^{**}$ can be obtained, for example, by solving the optimization problem (\ref{Eq: Main_OCP}) without safety constraints (\ref{Eq: OCPSafetyCosnt}), or a forward simulation of the vehicles with, for example, a simple feedback controller. 

The MIQP problem \eqref{Eq: MIQP_SQP} can be compactly written as 
\begin{subequations}
\label{Eq: MIQP}
\begin{align}
    &\underset{\mathcal{W}, b}{\textup{min}} \;\; \frac{1}{2}\begin{bmatrix}
\mathcal{W} \\
b
\end{bmatrix}^T \mathbf{H} \begin{bmatrix}
\mathcal{W} \\
b
\end{bmatrix} + \mathbf{J}^T\begin{bmatrix}
\mathcal{W} \\
b
\end{bmatrix} + \alpha \\
    &\textup{s.t.} \;\;\;  A_{\textup{eq}}\begin{bmatrix}
\mathcal{W}\\
b
\end{bmatrix} = b_{\textup{eq}}\\
    & \hspace{0.7cm} A_{\textup{ineq}}\begin{bmatrix}
\mathcal{W} \\
b
\end{bmatrix} \leq b_{\textup{ineq}}  \label{Eq: MIQP_INEQ},
\end{align}
\end{subequations}
\noindent where $\mathbf{J}$ now contains all the first order terms, $\alpha$ contains the linear terms and where the constraints \eqref{Eq: SQP_Const1}-\eqref{Eq: SQP_Const3} are grouped into equality constraints $A_{\textup{eq}}, b_{\textup{eq}}$ and inequality constraints $A_{\textup{ineq}}, b_{\textup{ineq}}$, respectively. The solution to the MIQP problem provides the ``approximately optimal" crossing orders $\hat{\mathcal{O}}^\mathcal{I}$, $\hat{\mathcal{O}}^\mathcal{M}$ that is obtained from the values of the integer variables $b$. 
\subsection{Fixed-order NLP}

Removing all other integer solutions than that found by the heuristic, \eqref{Eq: Main_OCP} is reduced to an NLP. Obtaining the optimal state and control trajectories is thus found through solving the \textit{fixed-order coordination problem}
\begin{subequations}
\label{Eq: FixedOrderNLP}
\begin{align}
    \underset{x_{i,k},u_{i,k}}{\textup{min}} \;\; &\sum_{i = 1}^{N_a}J_i\left ( x_{i,k},u_{i,k} \right )  \\
    \textup{s.t} \;\;\; & (\ref{Eq: InitialGuess})-(\ref{Eq: OCPSafetyCosnt}), \; \forall i, \forall k \\
    &\mathcal{O}^\mathcal{I} = \hat{\mathcal{O}}^\mathcal{I}, \;\; \mathcal{O}^\mathcal{M}=\hat{\mathcal{O}}^\mathcal{M}
\end{align}
\end{subequations}
The two stage approximation approach is summarized in Algorithm \ref{Al: ApproximationAlgrthm}.
\begin{algorithm}
\caption{Two stage approximation algorithm}
\begin{flushleft}
        \textbf{Input:} $N_a,\mathcal{I}, \mathcal{Q}_r, \mathcal{M}, \mathcal{Z}_w$, vehicle paths\\
        \textbf{Output:} $\mathcal{X}^{*}, \; \mathcal{U}^{*}$
\end{flushleft}
\begin{algorithmic}[1] 
\State $\forall i$: Obtain a solution guess $w_i^{**}$ by, e.g., solving NLP (\ref{Eq: Main_OCP}) w/o safety constraints (\ref{Eq: IntersectionConst}). \label{Al: ApproximationAlgrthm}
\State Calculate and form the approximation terms $\mathbf{H}, \mathbf{J}, \alpha$.
\State Solve the MIQP (\ref{Eq: MIQP}) to get ``approximately optimal" crossing orders $\hat{\mathcal{O}}^\mathcal{I}$, $\hat{\mathcal{O}}^\mathcal{M}$.
\State Solve the fixed-order NLP (\ref{Eq: FixedOrderNLP}) using $\hat{\mathcal{O}}^\mathcal{I}$, $\hat{\mathcal{O}}^\mathcal{M}$ to obtain $\mathcal{X}^{*}, \; \mathcal{U}^{*}$. 
\end{algorithmic}
\end{algorithm}

\section{Simulation results}\label{Ch:Simulation results}
In this section, we present a simulation example showing the operation and performance of the coordination algorithm.
\subsection{Simulation setup}
The vehicles are modelled as a triple integrator $\dddot{x} = u$, whereby the spatial model is: 
\begin{equation} 
    \begin{bmatrix}
    \tder{t_i}{p_i} \\
    \tder{v_i}{p_i} \\
    \tder{a_i}{p_i}
    \end{bmatrix} = \begin{matrix}
\frac{1}{v_i}\\ 
\frac{a_i}{v_i}\\ 
\frac{u_i}{v_i}
\end{matrix},
\end{equation}
\noindent where $a_i$ is the acceleration and $u_i$ is the jerk $j_i$. The problem is transcribed using multiple shooting (with $N$ shooting points) and an Explicit Runge-Kutta-4 (ERK4) integrator.

The state constraints are chosen as bounds on the speed and longitudinal acceleration $(\underline{v}_i \leq v_{i,k} \leq \overline{v}_i, \; a_{i,k} \leq \overline{a}_{i,lon})$ in order to obey speed limits and physical constraints. The vehicles are expected to manoeuvre on curved roads, thus, it is necessary to limit the lateral forces to avoid vehicle stability problems, like roll over. As the one-dimensional model that is used in this paper does not account for lateral motion, the following constraint is enforced, as similarly proposed in \cite{b8},
\begin{equation}\label{Eq: StateConst2}
    \left ( \frac{a_{i,k}}{a_{i,lon}} \right )^2 + \left ( \frac{\kappa_i(s_{i,k})v_{i,k}}{\overline{a}_{i,lat}} \right )^2 \leq 1,
\end{equation}
\noindent where $\overline{a}_{i,lat}$ is the lateral acceleration limit and $\kappa_i(p_{i,k})$ is the road curvature, that is assumed to be available at every point along the path.  

\begin{figure}[h]
    \centering
    \includegraphics[width=0.45\textwidth,height=12.5cm]{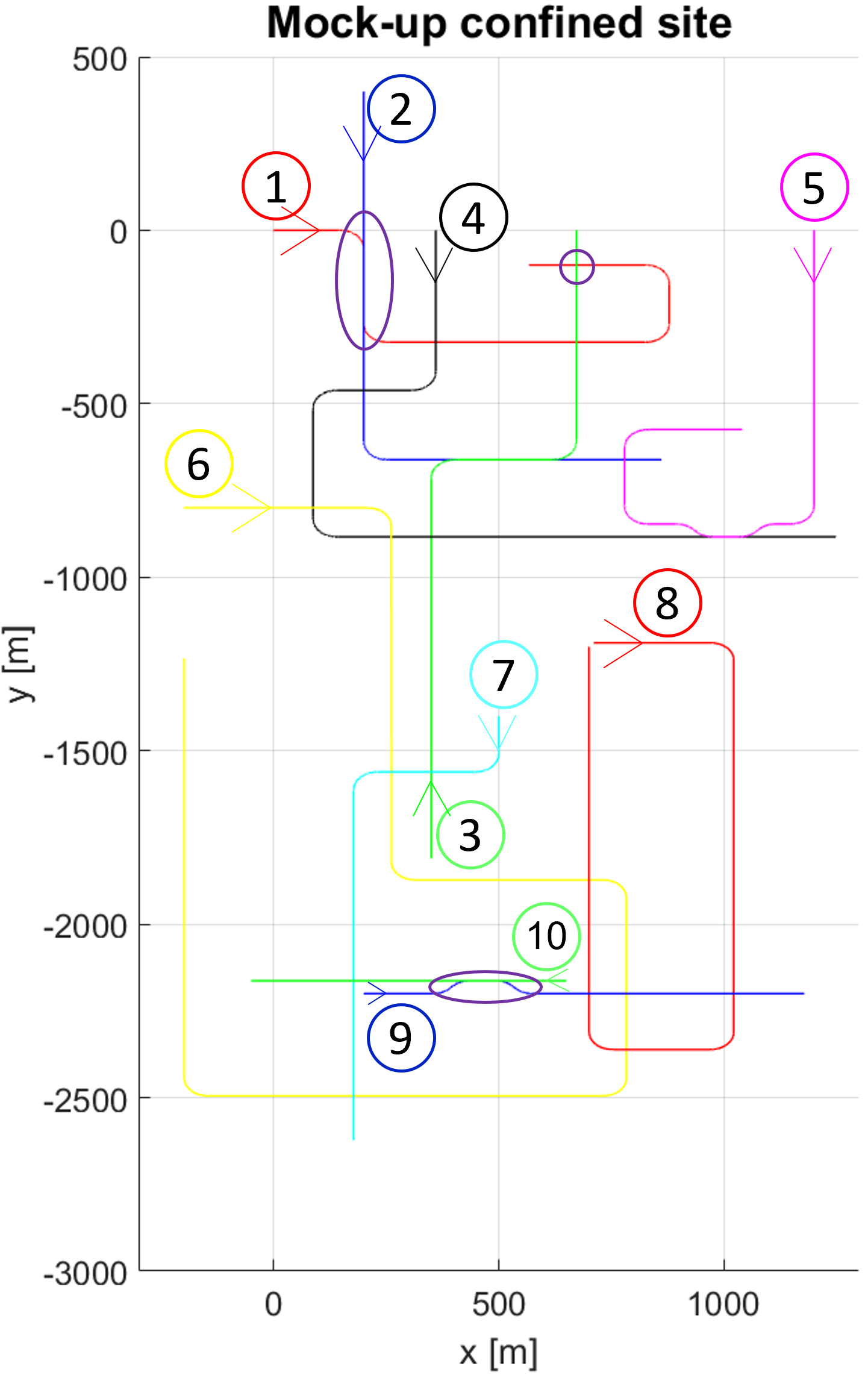}
    \centering
    \caption{Mock-up confined site area. The CZs that are investigated in the simulation scenario are circled over in this figure.}
    \label{Fig: SiteMap}
\end{figure}

While other objectives could be used, we consider a trade-off between minimization of the total travel time and the squares of the longitudinal acceleration and longitudinal jerk:
\begin{align}
&J(x_i(p_i),u_i(p_i)) = \\ \nonumber
&\int_{p_{i,0}}^{p_{i,N_i}}  \left ( P_ia_i(p_i)^2 + Q_ij_i(p_i)^2 \right )\frac{1}{v_i(p_i)}\mathrm{d}p_i + R_it_i(p_{i,N_i}),
\end{align}
\noindent where $P_i,Q_i,R_i$ are the appropriate weights. The inclusion of the acceleration and jerk in the objective can be interpreted as a drivability measure. The objective function is integrated with Forward Euler, leading to the ``discretized" objective function
\begin{align}
    &J_i\left ( x_{i,k},u_{i,k} \right ) = \\ \nonumber
    &\sum_{k=1}^{N}\left (\left ( P_ia_{i,k}^2 + Q_ij_{i,k}^2 \right ) \frac{\Delta p_{i,k}}{v_{i,k}}\right ) + R_it_{i,N},
\end{align}
\noindent with $\Delta p_{i,k} = p_{i,k+1}-p_{i,k}$.
We consider a confined site with layout shown in Figure \ref{Fig: SiteMap}. There are in total ten vehicles, two merge-split CZs, two narrow road CZs and sixteen intersection CZs. Every vehicle starts from an initial velocity of 50 $\textup{km/h}$ and vehicles 5 and 10 start from a nonzero initial time to ensure that a collision occurs if no coordinating action is taken. The respective initial times for those vehicles are $t_{0,5} = 56$ and $t_{0,10} = 11.3$ seconds. The chosen velocity and acceleration bounds are $\underline{v}_i = 3.6 \; \textup{km/h}$, $\overline{v}_i = 90 \; \textup{km/h}$, $\overline{a}_{i,lon} = 4 \; \textup{m}/\textup{s}^2$, $\overline{a}_{i,lat} = 2 \; \textup{m}/\textup{s}^2$. The weight coefficients for the objective are: $P_i = 1, \; Q_i = 1, \; R_i = 10$ and the number of shooting points is $N = 100$ for each vehicle. The intersection CZ is created with 5 meter margin ahead of and behind the collision point, where as in the merge-split and narrow road CZ the margin is 15 meters for both the entry and exit collision point. For the merge-split CZ it is desirable to keep at least 0.5s margin between the vehicles, i.e. $\Delta t= 0.5$ in \eqref{Eq:MergeSplit_ORConst}.

We utilize the CasADi \cite{b15} toolkit and  IPOPT \cite{b16}  to formulate and solve the optimization problem (\ref{Eq: FixedOrderNLP}) and use Gurobi for the MIQP (\ref{Eq: MIQP}).
\subsection{Discussion of results}
In the following, the results for this simulation scenario are presented. In particular, we demonstrate the uncoordinated and coordinated results for one merge-split zone, one intersection zone, and one narrow road, all circled over in Figure \ref{Fig: SiteMap}. The uncoordinated results are obtained as the individual optimum of each vehicle, i.e., the speed profiles are obtained by solving the optimization problem without any MUTEX zone constraints.
\begin{figure*}[h]
    \centering
    \begin{minipage}{0.48\textwidth}
        \centering
        \includegraphics[width=\textwidth,height=7cm]{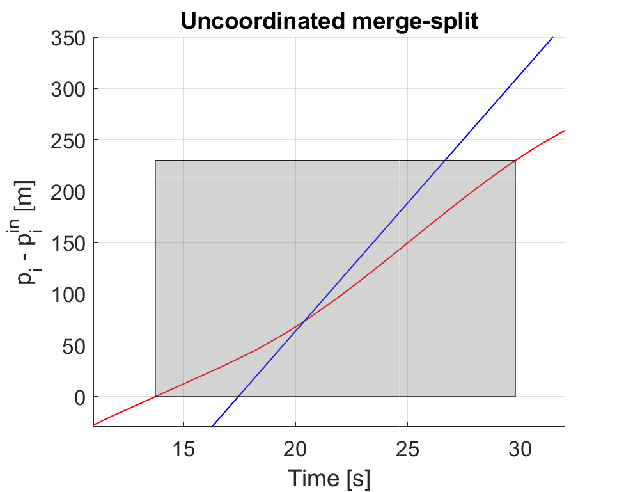} 
        \caption{Uncoordinated crossing for the merge-split CZ for the 1$^{\textup{st}}$ and 2$^{\textup{nd}}$ vehicle. The gray box is the occupancy time in the zone for the 1$^{\textup{st}}$ vehicle. While it is allowed for the 2$^{\textup{nd}}$ vehicle to enter the CZ whilst the 1$^{\textup{st}}$ vehicle is in it, an intersection between the vehicles translates to a collision between the vehicles.}
        \label{Fig: UncoordinatedMergeSplit}
    \end{minipage}\hfill
    \begin{minipage}{0.48\textwidth} 
        \centering
        \includegraphics[width=\textwidth,height=7cm]{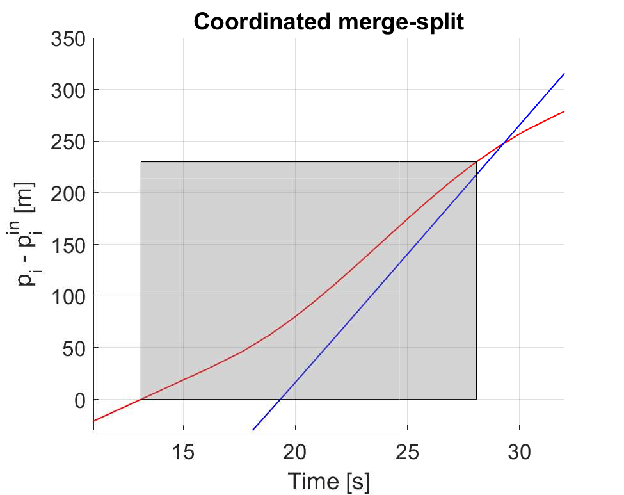}
        \caption{Coordinated crossing for the merge-split CZ for the 1$^{\textup{st}}$ and 2$^{\textup{nd}}$ vehicle. The gray box is the occupancy time in the zone for the 1$^{\textup{st}}$ vehicle. The 2$^{\textup{nd}}$ vehicle enters the CZ whilst the 1$^{\textup{st}}$ vehicle is in and keeps the desired time gap.}
        \label{Fig: CoordinatedMergeSplit}
    \end{minipage}
\end{figure*}

Figure \ref{Fig: UncoordinatedMergeSplit} and Figure \ref{Fig: CoordinatedMergeSplit} show the position vs. time trajectories of the two vehicles at the merge-split CZ in the uncoordinated and coordinated case, respectively. The trajectories are shifted such that position zero is the entry position of the CZ for both vehicles. The gray rectangle depicts the CZ and time spent inside the CZ for the 1$^{\textup{st}}$ vehicle. The interpretation of constraints \eqref{Eq:MergeSplit_ORConst} is that a collision occurs if the trajectories intersect while both vehicles are in the CZ. Figure \ref{Fig: UncoordinatedMergeSplit} shows that without coordination, the second vehicle would run into the first just after $t=20$ s (the trajectories intersect). In the coordinated case the vehicles are instead controlled to avoid a collision (the trajectories do not intersect) and keep at least the specified distance at all times. Note that an intersection between the trajectories outside the CZ is not relevant as they are no longer on a shared road. By allowing both vehicles to be in the CZ, and not blocking the whole zone for one vehicle, the throughput is increased.

\begin{figure}
    \centering
    \includegraphics[width=0.5\textwidth,height=9.5cm,keepaspectratio]{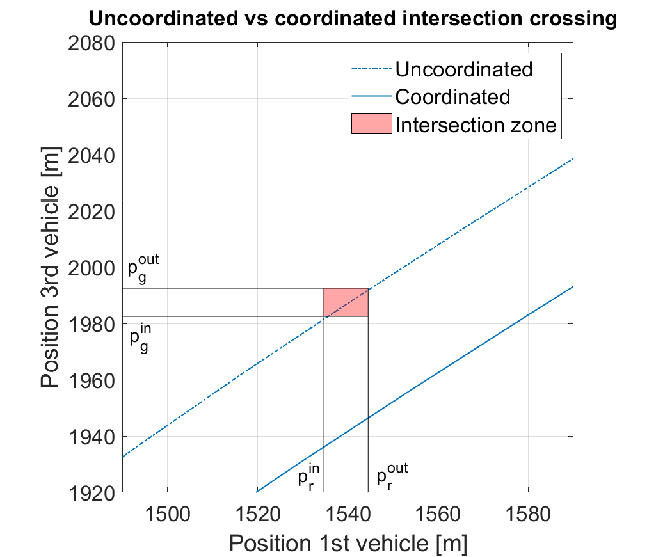}
    \caption{Uncoordinated and coordinated crossing of the intersection zone. An intersection of the depicted trajectory with the intersection zone is equivalent to both vehicles being in the CZ at the same time.}
    \label{Fig: IntersectionCrossing}
\end{figure}

In the scenarios where the paths of two (or more) vehicles intersect, it is necessary to only have one vehicle in the zone at any given time. The reasoning for this decision is that the CZ in this case occupies a relatively small patch of road and blocking off the whole zone for one vehicle will not result in a major loss of throughput. Figure \ref{Fig: IntersectionCrossing} shows the position-vs-position trajectory of the 1$^{st}$ and 3$^{rd}$ vehicle in the uncoordinated and coordinated case, respectively. The red area represents the CZ, but in this case, a collision occurs when two vehicles are inside the CZ at the same time. This is equivalent to all configurations where the trajectory is inside the gray area. As the figures show, a collision occurs in the uncoordinated case (the trajectory passes through the red area), the coordinated vehicles satisfy the collision constraints and avoid being in the CZ at the same time.  

\begin{figure}
    \centering
    \includegraphics[width=0.5\textwidth,height=9.5cm,keepaspectratio]{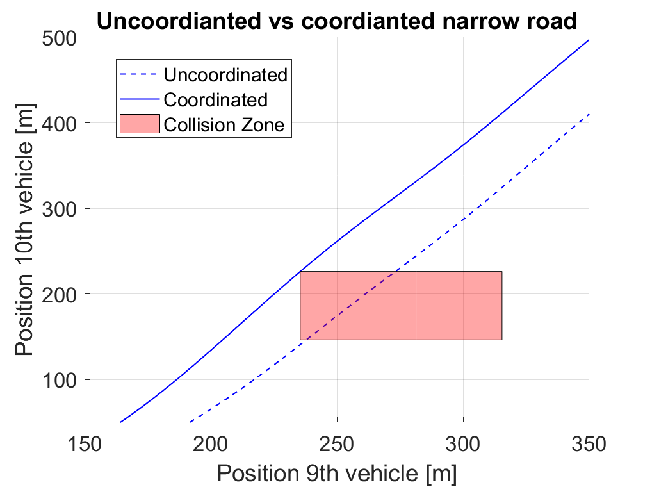}
    \caption{Uncoordinated and coordinated crossing of a narrow road collision zone. An intersection of the depicted trajectory with the CZ is equivalent to both vehicles being in the CZ at the same time.}
    \label{Fig: NarrowRoad}
\end{figure}

For the narrow road collision zones, the collision avoidance is defined in the same way as in the intersection zone since the road conditions allow for only one vehicle to occupy the zone at a time. Figure \ref{Fig: NarrowRoad} shows that both vehicles are inside the zone at the same time in the uncoordinated case (the trajectory passes through the red area). This is avoided in the coordinated case, although with a very small margin. The interpretation of the trajectory touching the corner of the red area is that one vehicle enters the narrow road just as the other exits.

\begin{figure}
    \centering
    \includegraphics[width=0.45\textwidth,height=26.5cm,keepaspectratio]{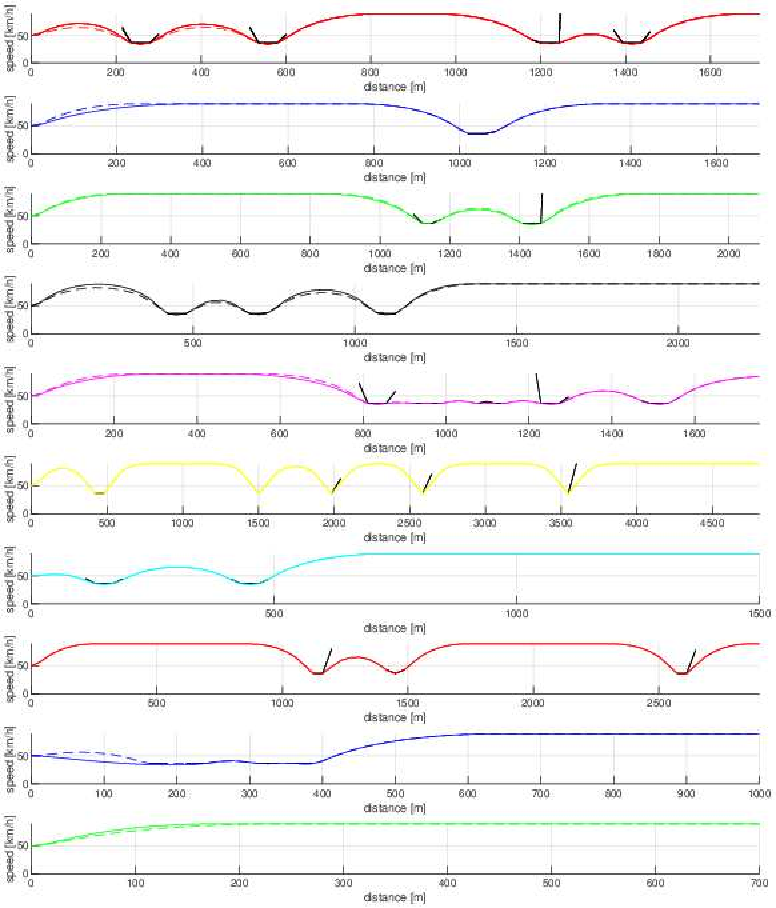}
    \caption{Speed profiles for the uncoordinated (dashed lines) and coordinated (solid lines) vehicles.}
    \label{fig:VehicleSpeed}
\end{figure}

Figure \ref{fig:VehicleSpeed} shows the speed profiles of the vehicles for the uncoordinated and coordinated case, where the black lines illustrate implicit bounds on speed through the acceleration constraint \eqref{Eq: StateConst2}. As the algorithm is aware of all CZs the vehicles encounters, it is able to avoid collisions with small speed changes well before the vehicle arrives at the CZ. With the objective function aiming to minimize the acceleration and jerk, it is noticeable that the speed profiles are smooth throughout the whole site and would not be challenging to follow in an actual application. It is worth to highlighting that several collisions occur in the uncoordinated case. When the vehicles are coordinated, all collisions are avoided. For sake of brevity, we refrain from depicting these results as the critical behaviour is similar to the cases discussed in detail above.

The simulation scenario is implemented in MATLAB on a 2.90GHz Intel Xeon computer with 32GB of RAM. The total solve computational time for the example is 1.973 seconds. To be precise, the MIQP requires 0.14 seconds and the main coordination optimal control problem requires 1.833 seconds.

\section{Conclusions and future work}\label{Ch:Conclusions}

In this paper, we have presented a heuristic that finds approximate solutions to the optimal coordination problem of automated vehicles at confined sites with multiple collision zones of different types. The approach optimizes the trajectories of the vehicles over their entire path while taking all collision zones into account. In future work, we intend to investigate improvement on the NLP problem computation time, closed-loop behaviour, and coupling the approach with a safety supervisor to guarantee satisfaction of the safety constraints.

\end{document}